\documentclass{amsart}
\usepackage{amssymb, amsmath}
\newcommand\field[1]{\mathbb{#1}}
\newcommand{\CC}{\field{C}}
\newcommand{\NN}{\field{N}}

\newcommand{\TT}{\field{T}}
\newcommand{\ZZ}{\field{Z}}
\newcommand{\Bb}{\mathcal B}

\newcommand\lsp{\operatorname{span}\nolimits}
\newcommand\clsp{\overline{\lsp}}

\theoremstyle{plain}
\newtheorem{theorem}{Theorem}[section]

\newtheorem{lemma}[theorem]{Lemma}
\newtheorem{prop}[theorem]{Proposition}
\theoremstyle{definition}
\newtheorem{dfn}{Definition}

\numberwithin{equation}{section}

\title[Simplicity of $k$-graph algebras]{Simplicity of $C^*$-algebras associated to higher-rank graphs}
\author{David I. Robertson}
\email{D.Robertson@newcastle.edu.au}
\address{\textnormal{Current address\footnotemark[1]:}
David Robertson\\
School of Mathematical and Physical Sciences\\
University of Newcastle\\
Callaghan \\
NSW 2308\\
AUSTRALIA}

\author{Aidan Sims}
\email{asims@uow.edu.au}
\address{\textnormal{Current address\footnotemark[1]:} Aidan Sims \\
School of Mathematics and Applied Statistics \\
Austin Keane building (15) \\
University of Wollongong \\
NSW 2522 \\
AUSTRALIA
}
\subjclass{Primary 46L05}
\keywords{Graph algebra, k-graph}
\thanks{The second author was supported by the Australian Research Council.}
\date{February 7, 2006}

\renewcommand{\theenumi}{\roman{enumi}}

\begin{document}
\newcommand{\itemref}[1]{\eqref{#1}}
\newcommand{\proofof}[1]{Proof of~#1}

\maketitle

\begin{abstract}
We prove that if $\Lambda$ is a row-finite $k$-graph with no
sources, then the associated $C^*$-algebra is simple if and
only if $\Lambda$ is cofinal and satisfies Kumjian and Pask's
aperiodicity condition, known as Condition~(A). We prove that
the aperiodicity condition is equivalent to a suitably modified
version of Robertson and Steger's original nonperiodicity
condition~(H3) which in particular involves only finite paths.
We also characterise both cofinality and aperiodicity of
$\Lambda$ in terms of ideals in $C^*(\Lambda)$.
\end{abstract}

\section{Introduction}

Consider a directed graph $E$ which is row-finite and has no
sinks in the sense that the set of outgoing edges from each
vertex is both finite and nonempty. As in~\cite{KPR}, we say
$E$ satisfies Condition~(L) if every cycle in $E$ has an exit,
and is cofinal if every vertex connects to every infinite path.
There is an elegant relationship between these conditions and
the ideal-structure of the graph algebra $C^*(E)$ (see
\cite{CBMSbk} for an overview). In particular:
\begin{enumerate}\renewcommand{\theenumi}{\arabic{enumi}}
\item\label{item:1-graph aperiodicity} every ideal of
    $C^*(E)$ contains at least one of the canonical
    generators of $C^*(E)$ if and only if $E$ satisfies
    Condition~(L) \cite[Theorem~3.7]{KPR}; and
\item\label{item:1-graph cofinality} $C^*(E)$ is simple if
    and only if $E$ satisfies Condition~(L) and is cofinal
    \cite[Proposition~5.1]{BPRS}.
\end{enumerate}

The $k$-graphs developed by Kumjian and Pask in \cite{KP} are a
generalisation of directed graphs designed to model the
higher-rank Cuntz-Krieger algebras of \cite{RobSt}. In
\cite{KP}, Kumjian and Pask identified a generalisation of
Condition~(L) for higher-rank graphs which they called the
\emph{aperiodicity condition} or Condition~(A), and showed that
this condition guarantees that every ideal of the $C^*$-algebra
contains at least one of the canonical generators. They also
identified a cofinality condition on higher-rank graphs which
together with the aperiodicity condition implies that the
associated $C^*$-algebra is simple. These two results
generalise the ``if'' directions of statements
\itemref{item:1-graph aperiodicity}~and~\itemref{item:1-graph
cofinality} of the previous paragraph. However, the
generalisations to $k$-graphs of the ``only if" directions of
\itemref{item:1-graph aperiodicity}~and~\itemref{item:1-graph
cofinality} above have not been established. Moreover, the
aperiodicity condition is phrased in terms of infinite paths,
and is difficult to verify in practise.

If we remove the hypotheses that a directed graph $E$ is
row-finite and has no sources, Condition~(L) and cofinality as
in \cite{KPR} still yield the same consequences for $C^*(E)$
\cite{DT2}, \cite{FLR}. For $k$-graphs, however, different
conditions (\cite[Condition~(B)]{RSY1} and
\cite[Condition~(C)]{RSY2}) have arisen as each hypothesis has
been removed. In the situation considered by Kumjian and Pask,
Conditions (B)~and~(C) are equivalent and imply aperiodicity
condition, but whether the reverse implication holds is an open
question.

A number of authors (see for example \cite{KP}, \cite{RSY1},
\cite{CBMSbk}) have pointed to the shortcomings of the
aperiodicity condition outlined above as significant open
problems. In this paper we resolve many of them for row-finite
$k$-graphs with no sources. In summary:
Theorem~\ref{thm:simple} shows that a row-finite $k$-graph with
no sources satisfies the aperiodicity condition and is cofinal
if and only if its $C^*$-algebra is simple. More specifically,
Proposition~\ref{prp:cofinal necessary} establishes
statement~\itemref{item:1-graph cofinality} above for
row-finite $k$-graphs with sources, and
Proposition~\ref{prp:aperiodicity necessary}
establishes~\itemref{item:1-graph aperiodicity} with
Condition~(L) replaced by the aperiodicity condition.
Additionally, Proposition~\ref{prp:all ideals g-i}
characterises the $k$-graphs for which every ideal of
$C^*(\Lambda)$ is gauge-invariant (c.f.
\cite[Section~6]{KPRR}). In Lemma~\ref{lem:equiv aperiodicity},
we identify a generalisation of Robertson and Steger's
nonperiodicity condition~(H3) to row-finite $k$-graphs with no
sources, and show that this condition, the aperiodicity
condition, and Condition~(B) are all equivalent. Significantly,
our generalisation of~(H3) involves only finite paths.
\vskip0.25em The $k$-graph versions of cofinality,
aperiodicity, the Cuntz-Krieger uniqueness theorem, and the
simplicity theorem for $k$-graphs used in this paper are all
due to Kumjian and Pask in~\cite{KP}. Nonetheless, we have
generally referenced~\cite{RSY1} throughout. This should not be
interpreted as a dismissal of the ground-breaking work of
Kumjian and Pask, which is undoubtedly the original and
definitive article. Our choice was made only
because~\cite{RSY1} is the earliest single paper containing
both the uniqueness theorems \emph{and} a description of the
gauge-invariant ideal structure for $k$-graph algebras, and
hence allowed us to restrict our referencing to a single paper.
\vskip0.25em \textbf{Acknowledgements.} We would like to thank
Trent Yeend for drawing our attention to the elementary proof
that Condition~(B) and the aperiodicity condition are
equivalent.

\section{Preliminaries}\label{sec:prelims}

In this section we gather the notation and conventions we need
regarding $k$-graphs. For a more detailed and rigorous
treatment of $k$-graphs, see \cite{KP}, \cite{RSY1}.

We regard $\NN^k$ as a monoid under addition, and denote its
generators by $e_1, \dots, e_k$. We write $n_i$ for the $i^{\rm
th}$ coordinate of $n \in \NN^k$. For $m,n \in \NN^k$, we say
$m\leq n$ if $m_i\leq n_i$ for each $i$. We write $m\vee n$ for
the coordinate-wise maximum of $m$ and~$n$. \vskip0.25em
\textbf{Higher-rank graphs.} Fix $k > 0$. We think of a
$k$-graph as a collection $\Lambda$ of paths endowed with a
\emph{degree} function $d : \Lambda \to \NN^k$ such that
concatenation of paths has the following \emph{factorisation
property}: if $d(\lambda)=m+n$, there are unique paths $\mu\in
d^{-1}(m)$ and $\nu\in d^{-1}(n)$ such that $\lambda=\mu\nu$.
For $n \in \NN^k$, we write $\Lambda^n$ for $d^{-1}(n)$. One
can make this notion rigorous using category-theory; see
\cite[Definition~1.1]{KP} for the formal definition.

The vertices of $\Lambda$ are the paths of degree $0$. For a
given path $\lambda$, the factorisation property ensures that
there are unique vertices, called the range and source of
$\lambda$ and denoted $r(\lambda)$ and $s(\lambda)$, such that
$r(\lambda) \lambda = \lambda = \lambda s(\lambda)$. For $v \in
\Lambda^0$ and $E\subseteq\Lambda$, we write $v E$ for $E \cap
r^{-1}(v)$ and $E v$ for $E \cap s^{-1}(v)$.

If $\lambda\in\Lambda$ with $d(\lambda)=l$, and $0\leq m\leq
n\leq l$, there exist unique paths $\lambda(0,m) \in
\Lambda^m$, $\lambda(m,n) \in \Lambda^{n-m}$ and $\lambda(n,l)
\in \Lambda^{l-n}$ such that $\lambda =
\lambda(0,m)\lambda(m,n)\lambda(n,l)$.

We say that a $k$-graph $\Lambda$ is \emph{row-finite} if
$v\Lambda^n$ is finite for all $v \in \Lambda^0$ and $n \in
\NN^k$, and we say $\Lambda$ has \emph{no sources} if
$v\Lambda^n$ is always nonempty.
\vskip0.25em \textbf{Infinite paths, the shift map, and
cofinality, and aperiodicity.} We denote by $\Omega_k$ the
$k$-graph with paths $\{(m,n)\in\NN^k\times\NN^k:m\leq n\}$,
$r(m,n)=(m,m)$, $s(m,n)=(n,n)$, $(m,n)(n,p)=(m,p)$, and
$d(m,n)=n-m$. For brevity, we generally write $n$ for the
vertex $(n,n)$ of $\Omega_k$.

Given $k\in\NN\backslash\{0\}$ and $k$-graphs $\Lambda$ and
$\Gamma$, a \emph{graph morphism} from $\Lambda$ to $\Gamma$ is
a function $x : \Lambda \to \Gamma$ which respects both
connectivity and degree. Given a $k$-graph $\Lambda$, an
\textit{infinite path} in $\Lambda$ is a graph morphism
$x:\Omega_k\to\Lambda$. We write $\Lambda^{\infty}$ for the
collection of all infinite paths in $\Lambda$. We denote $x(0)$
by $r(x)$, and call it the \emph{range} of $x$, and for $v \in
\Lambda^0$, we write $v\Lambda^\infty$ for the set $\{x \in
\Lambda^\infty : r(x) = v\}$. For $p\in\NN^k$, we write
$\sigma^p:\Lambda^{\infty}\to\Lambda^{\infty}$ for the shift
map determined by $\sigma^p(x)(n)=x(n+p)$.

We say $\Lambda$ is \emph{cofinal} if, for every $x \in
\Lambda^\infty$ and $v \in \Lambda^0$ there exists $n \in
\NN^k$ such that $v \Lambda x(n)$ is nonempty.

As in \cite{KP}, we say that $\Lambda$ satisfies the
\emph{aperiodicity condition} (or is \emph{aperiodic} if, for
every vertex $v \in \Lambda^0$ there is an infinite path $x \in
v\Lambda^\infty$ such that $\sigma^m(x) \not= \sigma^n(x)$ for
all $m \not= n \in \NN^k$. The aperiodicity condition is also
referred to as Condition~(A) in \cite{KP} and in other
$k$-graph literature. As in \cite{RSY1}, we say that $\Lambda$
satisfies \emph{Condition~(B)} if for every vertex $v \in
\Lambda^0$, there is an infinite path $x \in v\Lambda^\infty$
such that $\mu x \not= \nu x$ for all $\mu \not= \nu \in
\Lambda v$.
\vskip0.25em \textbf{$C^*$-algebras of $k$-graphs.} Let
$\Lambda$ be a row-finite $k$-graph with no sources. The
associated $C^*$-algebra $C^*(\Lambda)$ is the universal
$C^*$-algebra generated by partial isometries $\{s_\lambda :
\lambda \in \Lambda\}$ which satisfy the Cuntz-Krieger
relations:
\begin{enumerate}\renewcommand{\theenumi}{\arabic{enumi}}
\item\label{item:CK1} $\{s_v : v \in \Lambda^0\}$ are
    mutually orthogonal projections;
\item\label{item:CK2} $s_\mu s_\nu = s_{\mu\nu}$ when
    $s(\mu) = r(\nu)$;
\item\label{item:CK3} $s^*_\lambda s_\lambda =
    s_{s(\lambda)}$ for all $\lambda \in \Lambda$; and
\item\label{item:CK4} $s_v = \sum_{\lambda \in v\Lambda^n}
    s_\lambda s^*_\lambda$ for all $v \in \Lambda^0$ and $n
    \in \NN^k$.
\end{enumerate}
By \emph{universal} we mean that if $\{t_\lambda : \lambda \in
\Lambda\}$ is any collection of partial isometries in a
$C^*$-algebra $A$ which satisfy the Cuntz-Krieger
relations~\itemref{item:CK1}--\itemref{item:CK4} above, then
there is a homomorphism $\pi_t : C^*(\Lambda) \to A$ satisfying
$\pi_t(s_\lambda) = t_\lambda$ for all $\lambda \in \Lambda$.
Proposition~2.11 of \cite{KP} implies that the generators
$\{s_\lambda : \lambda \in \Lambda\}$ of $C^*(\Lambda)$ are all
nonzero.

\section{Results}\label{sec:results}

\begin{dfn}\label{dfn:no local p}
Let $\Lambda$ be a row-finite $k$-graph with no sources, and
let $v \in \Lambda^0$. We say that $\Lambda$ has \emph{local
periodicity at $v$} if there exist $m \not=n \in \NN^k$ such
that $\sigma^m(x) = \sigma^n(x)$ for all $x \in
v\Lambda^\infty$. We say that $\Lambda$ has \emph{no local
periodicity} if for each $v \in \Lambda^0$ and each $m \not= n
\in \NN^k$ there exists $x \in v\Lambda^\infty$ such that
$\sigma^m(x) \not= \sigma^n(x)$.
\end{dfn}

Note that the hypothesis that $\Lambda$ has no local
periodicity is only a slight weakening of the aperiodicity
condition: if $\Lambda$ satisfies the aperiodicity condition,
then for fixed $v \in \Lambda^0$ and distinct $m,n \in \NN^k$,
the path $x$ with range $v$ such that $\sigma^p(x) \not=
\sigma^q(x)$ whenever $p \not= q$ certainly satisfies
$\sigma^m(x) \not= \sigma^n(x)$. Indeed, we shall show in
Lemma~\ref{lem:equiv aperiodicity} that they are equivalent.

\begin{theorem}\label{thm:simple}
Let $\Lambda$ be a row-finite $k$-graph with no sources. Then
$C^*(\Lambda)$ is simple if and only if both of the following
conditions hold:
\begin{enumerate}
\item\label{item:thm cofinal} $\Lambda$ is cofinal; and
\item\label{item:thm aperiodic} $\Lambda$ has no local
    periodicity.
\end{enumerate}
\end{theorem}

We first show that the no local periodicity hypothesis, the
aperiodicity condition, and \cite[Condition~(B)]{RSY1} are all
equivalent to a version of \cite[(H3)]{RobSt} which involves
only finite paths of $\Lambda$.

\begin{lemma}\label{lem:equiv aperiodicity}
Let $\Lambda$ be a row-finite $k$-graph with no sources. Then
the following are equivalent.
\begin{enumerate}
\item\label{item:conditionA} $\Lambda$ satisfies the
    aperiodicity condition.
\item\label{item:conditionB} $\Lambda$ satisfies
    Condition~(B).
\item\label{item:inf-p np} $\Lambda$ has no local
    periodicity.
\item\label{item:fin-p np} For each vertex $v \in
    \Lambda^0$ and each pair $m \not= n \in \NN^k$ there is
    a path $\lambda = \lambda_{v,m,n} \in v\Lambda$ such
    that $d(\lambda) \ge m \vee n$  and
\[
\lambda\big(m, m + d(\lambda) - (m \vee n)\big) \not= \lambda\big(n, n+d(\lambda) - (m \vee n)\big).
\]
\end{enumerate}
\end{lemma}

Condition~\itemref{item:fin-p np} of Lemma~\ref{lem:equiv
aperiodicity} is a generalisation of the nonperiodicity
condition~(H3) of \cite{RobSt}. It looks complicated, but as a
formulation of aperiodicity in terms of finite paths, it is an
important outcome of the paper. Consequently we give a
pictorial explanation of the condition in
Appendix~\ref{app:fin-p np}.

\begin{proof}[\proofof{equivalence of \itemref{item:conditionA}, \itemref{item:inf-p np}~and~\itemref{item:fin-p np} in Lemma~\ref{lem:equiv aperiodicity}}]
\itemref{item:conditionA}$\implies$\itemref{item:inf-p np} was
    established in the remark immediately following
    Definition~\ref{dfn:no local p}.

\itemref{item:inf-p np}$\implies$\itemref{item:fin-p np}.
    Suppose $\Lambda$ has no local periodicity, and fix $v \in
    \Lambda^0$ and $m\not=n \in \NN^k$. Then there exists $x
    \in v\Lambda^\infty$ such that $\sigma^m(x) \not=
    \sigma^n(x)$. Hence $\sigma^m(x)(0,p) \not=
    \sigma^n(x)(0,p)$ for large enough $p \in \NN^k$. Let
    $\lambda_{v,m,n} := x(0, (m\vee n) + p)$. Then
\[
\lambda_{v,m,n}(m, m+p) = \sigma^m(x)(0,p) \not= \sigma^n(x)(0,p) = \lambda_{v,m,n}(n, n+p).
\]

\itemref{item:fin-p np}$\implies$\itemref{item:conditionA}.
    Suppose \itemref{item:fin-p np} holds, and fix $v \in
    \Lambda^0$. Let $(m^i, n^i)^\infty_{i=1}$ be a listing of
    $\{(m,n) \in \NN^k \times \NN^k : m \not= n\}$. Let
    $\lambda_1 := \lambda_{v, m^1, n^1}$ as
    in~\itemref{item:fin-p np}, and inductively let $\lambda_i
    := \lambda_{s(\lambda_{i-1}), m^i, n^i}$. Let $\eta_i :=
    \lambda_1 \lambda_2 \dots \lambda_i$ for each $i \ge 1$.
    Since $\bigvee_{i \in \NN} d(\eta_i) = \infty^k$, there is
    a unique infinite path $x \in \Lambda^\infty$ such that
    $x(0, d(\eta_i)) = \eta_i$ for all $i$
    (see~\cite[page~107]{RSY1}). Fix $m \not= n \in \NN^k$.
    Then $(m,n) = (m^i, n^i)$ for some $i$, and then
\begin{align*}
\sigma^{d(\eta_{i-1}) + m}(x)&(0, d(\lambda_i) - (m \vee n)) \\
&= \sigma^{d(\eta_{i-1})}(x)(m^i, m^i + d(\lambda_i) - (m^i \vee n^i)) \\
&= \lambda_i(m^i, m^i + d(\lambda_i) - (m^i \vee n^i))
\end{align*}
and likewise
\[
\sigma^{d(\eta_{i-1}) + n}(x)(0, d(\lambda_i) - (m \vee n)) = \lambda_i(n^i, n^i + d(\lambda_i) - (m^i \vee n^i))
\]
Hence $\sigma^m(x) \not= \sigma^n(x)$ by definition of
$\lambda_i$. Since $m,n$ were arbitrary, $x$ is aperiodic, and
since $v$ was arbitrary, it follows that $\Lambda$ satisfies
the aperiodicity condition.
\end{proof}

To prove that Condition~(B) is equivalent to the other
conditions in Lemma~\ref{lem:equiv aperiodicity}, we use a
technical lemma which we will use again in the proof of
Proposition~\ref{prp:aperiodicity necessary}.

\begin{lemma}\label{lem:technical}
Let $\Lambda$ be a row-finite $k$-graph with no sources.
Suppose that $\Lambda$ has local periodicity at $v$, and fix $m
\not= n \in \NN^k$ such that $\sigma^m(x) = \sigma^n(x)$ for
all $x \in v\Lambda^\infty$. Fix $\mu \in v\Lambda^m$ and
$\alpha \in s(\mu) \Lambda^{(m \vee n) - m}$ and let $\nu =
(\mu\alpha)(0,n)$. Then $\mu\alpha y = \nu\alpha y$ for all $y
\in s(\alpha)\Lambda^\infty$.
\end{lemma}
\begin{proof}
Fix $y \in s(\alpha) \Lambda^\infty$, and let $x := \mu\alpha
y$. Since $x \in v\Lambda^\infty$, $\sigma^m(x) = \sigma^n(x)$.
We have $\sigma^m(x) = \alpha y$ by definition. Since $\nu =
(\mu\alpha)(0,n)$, $x(0,n) = \nu$, so $x = \nu \sigma^n(x)$. As
$\sigma^n(x) = \sigma^m(x)$, it follows that $\sigma^n(x) =
\alpha y$, so $\mu\alpha y = x = \nu\sigma^n(x) = \nu\alpha y$.
\end{proof}

\begin{proof}[\proofof{equivalence of \itemref{item:conditionA}~and~\itemref{item:conditionB} in Lemma~\ref{lem:equiv aperiodicity}}]
Remark~4.4 of \cite{RSY1} shows that the aperiodicity condition
implies Condition~(B). Since \itemref{item:conditionA},
\itemref{item:inf-p np}~and~\itemref{item:fin-p np} are
equivalent, it now suffices to show that Condition~(B) implies
that $\Lambda$ has no local periodicity. We argue by
contrapositive. Suppose $\Lambda$ has local periodicity at $v$.
Let $m$, $n$, $\mu$, $\nu$ and $\alpha$ be as in
Lemma~\ref{lem:technical}. Since $d(\mu\alpha) = m + d(\alpha)
\not= n + d(\alpha) = d(\nu\alpha)$, we have $\mu\alpha \not=
\nu\alpha$. Since Lemma~\ref{lem:technical} implies that
$\mu\alpha y = \nu\alpha y$ for all $y \in
s(\alpha)\Lambda^\infty$, it follows that $\Lambda$ does not
satisfy Condition~(B).
\end{proof}

Our next steps are to describe the significance of cofinality
and of the aperiodicity condition in terms of ideals in
$C^*(\Lambda)$. The proof of Proposition~\ref{prp:cofinal
necessary} is not new (see for example
\cite[Proposition~5.1]{BPRS} and \cite[Proposition~4.8]{KP})
but the result has not to our knowledge been stated explicitly
before now.

For $z \in \TT^k$ and $n \in \ZZ^k$, we use the multi-index
notation $z^n$ for the product $z_1^{n_1} z_2^{n_2} \cdots
z_k^{n_k} \in \TT$. Recall from \cite[Section~4]{RSY1} that the
universal property of $C^*(\Lambda)$ supplies automorphisms
$\{\gamma_z : z \in \TT^k\}$ of $C^*(\Lambda)$ which satisfy
$\gamma_z(s_\lambda) = z^{d(\lambda)}s_\lambda$ for all
$\lambda \in \Lambda$ and $z \in \TT^k$. The map $z \mapsto
\gamma_z$ is a strongly continuous action of $\TT^k$ on
$C^*(\Lambda)$. The fixed point algebra $C^*(\Lambda)^\gamma$
is called the \emph{core} of $C^*(\Lambda)$ and is equal to
$\clsp\{s_\mu s^*_\nu : d(\mu) = d(\nu)\}$.

\begin{prop}\label{prp:cofinal necessary}
Let $\Lambda$ be a row-finite $k$-graph with no sources. The
following are equivalent:
\begin{enumerate}
\item\label{item:Lambda cofinal} $\Lambda$ is cofinal.
\item\label{item:vert proj gives everything} If $I$ is an
    ideal of $C^*(\Lambda)$ and $s_v \in I$ for some $v \in
    \Lambda^0$, then $I = C^*(\Lambda)$.
\item\label{item:intersects core gives everything} If $I$
    is an ideal of $C^*(\Lambda)$ and $I \cap
    C^*(\Lambda)^\gamma \not= \{0\}$,  then $I =
    C^*(\Lambda)$.
\end{enumerate}
\end{prop}

To prove the proposition, we need to recall some terminology
from \cite[Section~5]{RSY1}. We say that $H \subset \Lambda^0$
is \emph{hereditary} if $r(\lambda) \in H$ implies
$s(\lambda)\in H$ for all $\lambda \in \Lambda$, and that $H$
is \emph{saturated} if we have $v \in H$ whenever there exists
$n \in \NN^k$ such that $s(\lambda) \in H$ for all $\lambda \in
v\Lambda^n$.

\begin{proof}
We first show that \itemref{item:Lambda
cofinal}~and~\itemref{item:vert proj gives everything} are
equivalent, and then that \itemref{item:vert proj gives
everything}~and~\itemref{item:intersects core gives everything}
are equivalent.

\itemref{item:Lambda cofinal}$\implies$\itemref{item:vert proj
    gives everything}. Suppose that $\Lambda$ is cofinal. An
    argument formally identical to the second paragraph of
    \cite[Proposition~5.1]{BPRS} shows that the only nonempty
    saturated hereditary subset of $\Lambda^0$ is $\Lambda^0$
    itself. Theorem~5.2 of \cite{RSY1} then implies that the
    only ideal of $C^*(\Lambda)$ which contains a vertex
    projection is $C^*(\Lambda)$ itself.

\itemref{item:vert proj gives
    everything}$\implies$\itemref{item:Lambda cofinal}. We
    argue by contrapositive. Suppose that $\Lambda$ is not
    cofinal. We must construct an ideal $I$ of $C^*(\Lambda)$
    such that $I \not= C^*(\Lambda)$, but $s_v \in I$ for some
    $v \in \Lambda^0$.

Since $\Lambda$ is not cofinal,  there exists a vertex $v_0 \in
\Lambda^0$ and an infinite path $x \in \Lambda^\infty$ such
that $v_0 \Lambda x(n) = \emptyset$ for all $n \in \NN^k$. Let
\[
H_x := \{w \in \Lambda^0 : w\Lambda x(n) = \emptyset \text{ for all } n \in \NN^k\}.
\]
Then $\emptyset \subsetneq H_x \subsetneq \Lambda^0$ because
$v_0 \in H_x$ but $x(0) \not\in H_x$. We claim that $H_x$ is
saturated and hereditary.

To see that $H_x$ is hereditary, fix $u \in H_x$ and $v \in
\Lambda^0$ with $u\Lambda v \not= \emptyset$, say $\lambda \in
u\Lambda v$. If we suppose for contradiction that $v \not\in
H_x$, then there exists $\alpha \in v\Lambda x(n)$ for some $n
\in \NN^k$ and it follows that $\lambda\alpha \in u\Lambda
x(n)$ contradicting $u \in H_x$.

To see that $H_x$ is saturated, fix $v \in \Lambda^0$ and $m
\in \NN^k$ such that $s(\lambda) \in H_x$ for all $\lambda \in
v\Lambda^m$. Suppose for contradiction that $v \not\in H_x$.
Then there exists $\alpha \in v\Lambda x(n)$ for some $n \in
\NN^k$. Let $\alpha' = x(n, n + m)$. We have $s(\alpha') =
x(n+m)$ and $r(\alpha') = s(\alpha)$. By choice,
$d(\alpha\alpha') \ge m$ so we may rewrite $\alpha\alpha' =
\mu\nu$ where $\mu \in v\Lambda^m$. By choice of $m$, we have
$s(\mu) \in H_x$, and since $H_x$ is hereditary, it follows
that $s(\nu) \in H_x$. But $s(\nu) = s(\alpha') = x(n+m)$,
contradicting the definition of $H_x$.

Since $\emptyset \subsetneq H_x \subsetneq \Lambda^0$,
Theorem~5.2 of \cite{RSY1} implies that $I_{H_x}$ is a
nontrivial gauge-invariant ideal of $C^*(\Lambda)$.

\itemref{item:vert proj gives
    everything}$\implies$\itemref{item:intersects core gives
    everything}. Recall from \cite[Section~4]{RSY1} that there
    is an isomorphism $\pi$ from $C^*(\Lambda)^\gamma$ to
    $\varinjlim_{n \in \NN^k} \bigoplus_{v \in \Lambda^0}
    M_{\Lambda^n v}(\CC)$ which takes $s_\lambda s^*_\lambda$
    to the diagonal matrix unit $\theta_{\lambda,\lambda}$ in
    $M_{\Lambda^{d(\lambda)} s(\lambda)}(\CC)$. Let $I$ be an
    ideal of $C^*(\Lambda)$ which intersects
    $C^*(\Lambda)^\gamma$ nontrivially. Then $I$ must intersect
    $\pi^{-1}(\bigoplus_{v \in \Lambda^0} M_{\Lambda^n
    v}(\CC))$ for some $n$, hence must intersect one of the
    summands $\pi^{-1}(M_{\Lambda^n v})$. As
    $\pi^{-1}(M_{\Lambda^n v})$ is simple, it follows that $I
    \cap \pi^{-1}(M_{\Lambda^n v}) = \pi^{-1}(M_{\Lambda^n
    v})$, so $I$ contains $s_\lambda s^*_\lambda$ for some
    $\lambda \in \Lambda^n v$. Hence $s_v =
    s_\lambda^*(s_\lambda s^*_\lambda) s_\lambda \in I$, and
    \itemref{item:vert proj gives everything} implies that $I =
    C^*(\Lambda)$.

\itemref{item:intersects core gives
    everything}$\implies$\itemref{item:vert proj gives
    everything}. Trivial.
\end{proof}

Let $\{\xi_x : x \in \Lambda^\infty\}$ denote the usual basis
for $\ell^2(\Lambda^\infty)$. As in \cite[Theorem~3.15]{RSY1},
there is a family $\{S_\eta : \eta \in \Lambda\} \subset
\Bb(\ell^2(\Lambda^\infty))$ satisfying the Cuntz-Krieger
relations such that
\begin{equation}\label{eq:S_eta formula}
S_\eta \xi_x = \begin{cases}
\xi_{\eta x} &\text{ if $r(x) = s(\eta)$}\\
0 &\text{ otherwise,}
\end{cases}
\text{\quad and\quad}
S^*_\eta \xi_y = \begin{cases}
\xi_{\sigma^{d(\eta)}(y)} &\text{ if $y(0,d(\eta)) = \eta$} \\
0 &\text{ otherwise.}
\end{cases}
\end{equation}
The universal property of $C^*(\Lambda)$ gives a homomorphism
$\pi_S : C^*(\Lambda) \to \Bb(\ell^2(\Lambda^\infty))$
satisfying $\pi_S(s_\eta) = S_\eta$ for all $\eta \in \Lambda$.
We call $\pi_S$ the infinite path representation.

\begin{prop}\label{prp:aperiodicity necessary}
Let $\Lambda$ be a row-finite $k$-graph with no sources. The
following are equivalent:
\begin{enumerate}
\item\label{item:aperiodic} $\Lambda$ has no local
    periodicity.
\item\label{item:every ideal has vert proj} Every nonzero
    ideal of $C^*(\Lambda)$ contains a vertex projection.
\item\label{item:inf path repn faithful} The infinite path
    representation $\pi_S$ is faithful.
\end{enumerate}
\end{prop}
\begin{proof}
\itemref{item:aperiodic}$\implies$\itemref{item:every ideal has
    vert proj}. Suppose that for each $v \in \Lambda^0$ and
    each pair $m \not= n \in \NN^k$ there exists $x \in
    v\Lambda^\infty$ such that $\sigma^m(x) \not= \sigma^n(x)$.
    Then Lemma~\ref{lem:equiv aperiodicity} implies that
    $\Lambda$ satisfies Condition~(B). The Cuntz-Krieger
    uniqueness theorem \cite[Theorem~4.3]{RSY1} therefore
    implies that every ideal of $C^*(\Lambda)$ contains a
    vertex projection.

\itemref{item:every ideal has vert
    proj}$\implies$\itemref{item:inf path repn faithful}. For
    $v \in \Lambda^0$, $\pi_S(s_v) = S_v$ is the projection
    onto $\clsp\{\xi_x : x \in v\Lambda^\infty\}$ and so is
    nonzero. So $\ker(\pi_S)$ contains no vertex projection and
    is trivial by~\itemref{item:every ideal has vert proj}.

\itemref{item:inf path repn
    faithful}$\implies$\itemref{item:aperiodic}. We argue by
    contrapositive. Suppose that $\Lambda$ has local
    periodicity at $v \in \Lambda^0$. By
    Lemma~\ref{lem:technical} there exist $m \not= n \in
    \NN^k$, $\mu \in v\Lambda^m$, $\nu \in v\Lambda^n s(\mu)$
    and $\alpha \in s(\mu)\Lambda$ such that $\mu\alpha y =
    \nu\alpha y$ for all $y \in s(\alpha)\Lambda^\infty$.

We will show that $a := s_{\mu\alpha} s^*_{\mu\alpha} -
s_{\nu\alpha} s^*_{\mu\alpha}$ belongs to
$\ker(\pi_S)\setminus\{0\}$.

We begin by showing that $a$ is nonzero. The gauge action
$\gamma$ of $\TT^k$ on $C^*(\Lambda)$ satisfies
$\gamma_z(s_{\mu\alpha} s^*_{\mu\alpha}) = s_{\mu\alpha}
s^*_{\mu\alpha}$ and $\gamma_z(s_{\nu\alpha} s^*_{\mu\alpha}) =
z^{n-m} s_{\nu\alpha} s^*_{\mu\alpha}$. Fix $\omega \in \TT^k$
such that $\omega^{n-m} = -1$. Suppose for contradiction that
$a = 0$. Then
\begin{equation}\label{eq:generator zero}
0 = a + \gamma_\omega(a) = s_{\mu\alpha} s^*_{\mu\alpha} - s_{\nu\alpha} s^*_{\mu\alpha} + s_{\mu\alpha} s^*_{\mu\alpha} + s_{\nu\alpha} s^*_{\mu\alpha} = 2s_{\mu\alpha} s^*_{\mu\alpha},
\end{equation}
contradicting \cite[Theorem~3.15]{RSY1} which shows that
$s_\lambda \not= 0$ for all $\lambda \in \Lambda$.

To see that $\pi_S(a) = 0$, we fix $x \in \Lambda^\infty$ and
show that $\pi_S(a) \xi_x = 0$. By~\eqref{eq:S_eta formula},
\[
S_{\mu\alpha} S^*_{\mu\alpha} \xi_x =
\begin{cases}
\xi_x &\text{ if $x(0, m \vee n) = \mu\alpha$} \\
0 &\text{ otherwise.}
\end{cases}
\]
We consider two cases: either $x(0, m \vee n) \not =
\mu\alpha$, or $x = \mu\alpha y$ for some $y \in
\Lambda^\infty$. In the first case, we have $\pi_S(a) \xi_x =
0$ by~\eqref{eq:S_eta formula}. In the second case, our choice
of $\mu$, $\nu$, $\alpha$ ensures that $x = \mu\alpha y =
\nu\alpha y$. Equation~\eqref{eq:S_eta formula} therefore
implies that $S_{\nu\alpha} S^*_{\mu\alpha} \xi_x =
S_{\nu\alpha} \xi_y = \xi_x$. Hence $\pi_S(a)\xi_x = 0$.

As $x$ was arbitrary, $\pi_S(a)$ annihilates all basis elements
of $\ell^2(\Lambda^\infty)$, so is equal to zero. Since
$\pi_S(s_v) = S_v \not= 0$ for each $v \in \Lambda^0$,
$\ker(\pi_S)$ is an ideal of $C^*(\Lambda)$ which contains no
vertex projection.
\end{proof}

\begin{proof}[\proofof{Theorem~\ref{thm:simple}}]
If $C^*(\Lambda)$ is simple, then Proposition~\ref{prp:cofinal
necessary} implies that $\Lambda$ is cofinal, and
Proposition~\ref{prp:aperiodicity necessary} implies that
$\Lambda$ has no local periodicity. Conversely, if $\Lambda$ is
cofinal and has no local periodicity and $I$ is an ideal of
$C^*(\Lambda)$, then Proposition~\ref{prp:aperiodicity
necessary} implies $s_v \in I$ for some $v \in \Lambda^0$, and
then Proposition~\ref{prp:cofinal necessary} implies that $I =
C^*(\Lambda)$.
\end{proof}

Recall from \cite[Section~5]{RSY1} that if $H \subset
\Lambda^0$ is hereditary, then $\Lambda \setminus \Lambda H :=
\{\lambda \in \Lambda : s(\lambda) \not\in H\}$ is itself a
locally convex row-finite $k$-graph. It is easy to check that
if $H$ is also saturated, then $\Lambda \setminus \Lambda H$
also has no sources. As in \cite{RSY1}, given a saturated
hereditary $H \subset \Lambda^0$, we denote by $I_H$ the ideal
generated by $\{s_v : v \in H\}$; and given an ideal $I \subset
C^*(\Lambda)$, we denote by $H_I$ the collection $\{v \in
\Lambda^0 : s_v \in I\}$.

\begin{prop} \label{prp:all ideals g-i}
Let $\Lambda$ be a row-finite $k$-graph with no sources. Then
the following are equivalent:
\begin{enumerate}
\item\label{item:all ideals g-i} Every ideal of
    $C^*(\Lambda)$ is gauge-invariant.
\item\label{item:higher-rank condK} For every saturated
    hereditary $H \subset \Lambda^0$, $\Lambda \setminus
    \Lambda H$ has no local periodicity.
\end{enumerate}
\end{prop}
\begin{proof}
\itemref{item:higher-rank condK}$\implies$\itemref{item:all
    ideals g-i}. Lemma~\ref{lem:equiv aperiodicity} shows that
    each $\Lambda \setminus \Lambda H$ satisfies Condition~(B).
    Hence Theorem~5.3 of \cite{RSY1} implies that every ideal
    of $C^*(\Lambda)$ is gauge invariant.

\itemref{item:all ideals
    g-i}$\implies$\itemref{item:higher-rank condK}. We argue by
    contrapositive. Suppose that there is a saturated
    hereditary subset $H$ of $\Lambda$ such that $\Lambda
    \setminus \Lambda H$ has local periodicity at $v$, say. Let
    $\{t_\lambda : \lambda \in \Lambda \setminus \Lambda H\}$
    denote the universal generating Cuntz-Krieger family for
    $C^*(\Lambda \setminus \Lambda H)$. Theorem~5.2 of
    \cite{RSY1} shows that there is an isomorphism $\phi$ of
    $C^*(\Lambda)/I_H$ onto $C^*(\Lambda \setminus \Lambda H)$
    satisfying $\phi(s_\lambda + I_H) = t_\lambda$ for all
    $\lambda \in \Lambda \setminus \Lambda H$. Let $q_{I_H}$
    denote the quotient map from $C^*(\Lambda)$ to
    $C^*(\Lambda)/I_H$.

The argument of Proposition~\ref{prp:aperiodicity necessary}
gives a nonzero element $a = t_{\mu\alpha} t^*_{\mu\alpha} -
t_{\nu\alpha} t^*_{\mu\alpha}$ of $C^*(\Lambda \setminus
\Lambda H)$ which satisfies $\pi_T(a) = 0$ where $\pi_T$ is the
infinite-path representation of $C^*(\Lambda \setminus \Lambda
H)$. Let $b := s_{\mu\alpha} s^*_{\mu\alpha} - s_{\nu\alpha}
s^*_{\mu\alpha} \in C^*(\Lambda)$. By definition of $b$, we
have $\phi \circ q_{I_H}(b) = a \not= 0$, but $\pi_T \circ \phi
\circ q_{I_H}(b) = 0$. Since $\phi$ is an isomorphism, the
kernel of $\phi \circ q_{I_H}$ is precisely $I_H$. Theorem~5.2
of \cite{RSY1} implies that $H_{I_H} = H$. Since the kernel of
$\pi_T$ contains no of the vertex projections of $C^*(\Lambda
\setminus \Lambda H)$, the ideal $J = \ker(\pi_T \circ \phi
\circ q_{I_H})$ also satisfies $H_J = H$. Now $J \not= I_H$
because $b \in J \setminus I_H$. Theorem~5.2 of \cite{RSY1}
implies that $I \mapsto H_I$ is a bijection between
gauge-invariant ideals of $C^*(\Lambda)$ and saturated
hereditary subsets of $\Lambda^0$ with inverse $H \mapsto I_H$.
Since $J \not= I = I_{H_J}$, it follows that $J$ is a
nontrivial ideal of $C^*(\Lambda)$ which is not
gauge-invariant.
\end{proof}

\appendix
\section{Finite paths and aperiodicity}\label{app:fin-p np}
\setlength{\unitlength}{1em} Condition~\itemref{item:fin-p np}
of Lemma~\ref{lem:equiv aperiodicity} insists that for each $v
\in \Lambda^0$ and each $m \not= n \in \NN^k$ there exists a
path $\lambda \in v\Lambda$ such that $d(\lambda) \ge m \vee n$
and
\[
\lambda(m, m + d(\lambda) - (m \vee n)) \not= \lambda(n, n + d(\lambda) - (m \vee n)).
\]
In this appendix we attempt to provide some intuition for what
this condition says.

Fix a vertex $v$ in a $k$-graph $\Lambda$, and a pair $m \not= n \in \NN^k$. In Figure~\ref{fig:2-dim illustration}, a path $\lambda$ in $v\Lambda^{(n \vee m) + l}$ is represented by the large rectangle. %
\font\fourrm=cmr4
\begin{figure}[ht]
\[
\begin{picture}(13,9)
\put(0.7,0){$v$}
\put(1,9.3){\vector(0,-1){2.95}} \put(1,6.35){\vector(0,-1){4}} \put(1,2.35){\vector(0,-1){1.9}}
\put(13,0.3){\vector(-1,0){1.9}} \put(11.1,0.3){\vector(-1,0){5}} \put(6.1,0.3){\vector(-1,0){4.95}}
\put(13,9){\vector(0,-1){2.65}} \put(13,6.35){\vector(0,-1){0.7}} \put(13,5.15){\vector(0,-1){2.8}} \put(13,2.35){\vector(0,-1){2.05}}
\put(12.65,5.25){$m\!+\!l$}
\put(12.4,9.3){\vector(-1,0){1.25}} \put(11.15,9.3){\vector(-1,0){2.3}} \put(7.2,9.3){\vector(-1,0){1.1}} \put(6.1,9.3){\vector(-1,0){5.1}}
\put(12.5, 9.15){$(m\!\vee\!n)\!+\!l$}
\put(7.3,9.15){$n\!+\!l$}
\put(10.8,2.1){$m$}
\put(10.65,2.35){\vector(-1,0){9.65}}
\put(11.1,1.95){\vector(0,-1){1.65}}
\put(5.9,6.1){$n$}
\put(5.8,6.35){\vector(-1,0){4.8}}
\put(6.1,5.95){\vector(0,-1){5.65}}
\put(11.1,9.3){\vector(0,-1){2.7}}
\put(11.1,6){\vector(0,-1){3.4}}
\put(9.95, 6.35){\vector(-1,0){3.5}}
\put(10.1, 6.1){$m\!\vee\!n$}
\put(13, 6.35){\vector(-1,0){0.95}}
\put(6.1,9.3){\vector(0,-1){2.7}}
\put(8.1,9){\vector(0,-1){2.65}}
\put(6.85, 7.55){$\alpha$}
\put(13,2.35){\vector(-1,0){1.5}}
\put(12.55,5.4){\vector(-1,0){1.45}}
\put(11.75,3.55){$\beta$}
\multiput(11.125,2.425)(0,0.25){12}{
   \multiput(0,0)(0.25,0){8}{\fourrm.}
}
\multiput(11.25,2.55)(0,0.25){11}{
   \multiput(0,0)(0.25,0){7}{\fourrm.}
}
\multiput(6.125,6.425)(0,0.25){11}{
   \multiput(0,0)(0.25,0){8}{\fourrm.}
}
\multiput(6.25,6.55)(0,0.25){11}{
   \multiput(0,0)(0.25,0){7}{\fourrm.}
}
\end{picture}
\]
\caption{Condition~\itemref{item:fin-p np}, Lemma~\ref{lem:equiv aperiodicity} when $m, n$ are not comparable.}\label{fig:2-dim illustration}
\end{figure}%
Regarded as a scale diagram, Figure~\ref{fig:2-dim
illustration} illustrates the configuration $k = 2$, $m =
(10,2)$, $n = (5,6)$. However, we can use this picture to
represent the $k$-dimensional situation by using horizontal
distance to represent the directions in which $m$ is bigger
than $n$ and vertical distance to represent the directions in
which $n$ is bigger than~$m$.

If $\Lambda$ satisfies Condition~\itemref{item:fin-p np} of
Lemma~\ref{lem:equiv aperiodicity}, then there exists a path
$\lambda$ as in the diagram whose degree $(m \vee n) + l$ is
greater than both $m$ and $n$ and whose segment from $m$ to $m
+ l$ is distinct from the segment from $n$ to $n + l$. (In the
picture, $l = (2,3)$, but more generally it is the difference,
represented by the top-right rectangle, between the least upper
bound of $m$ and $n$ and the degree of $\lambda$.)

The factorisation property ensures that for each path from top
right to bottom left in the picture there is a unique
factorisation of $\lambda$ into segments of the corresponding
degrees. Condition~\itemref{item:fin-p np} of
Lemma~\ref{lem:equiv aperiodicity} insists that there exists
$\lambda$ as shown in Figure~\ref{fig:2-dim illustration} for
which the shaded segments $\alpha$ and $\beta$ are distinct.

\vskip2em \footnotetext[1]{Both authors were affiliated with
the University of Newcastle, Australia at the time of
publication.}
\end{document}